\documentclass[11pt,oneside,reqno]{article}
\usepackage{amsmath,amssymb,amsthm,mathrsfs,bbm}
\usepackage[round]{natbib}
\theoremstyle{plain}
\newtheorem{theorem}{Theorem}

\newtheorem{corollary}{Corollary}
\newtheorem{proposition}{Proposition}
\newcommand{\rnum}{\mathbb{R}}
\newcommand{\nnum}{\mathbb{N}}
\newcommand{\trace}{\mathop{\mathrm{tr}}}
\newcommand{\veps}{\varepsilon}
\newcommand{\dconv}{\stackrel{\mathcal{D}}{\rightarrow}}
\newcommand{\follows}{\,\Rightarrow\,}
\renewcommand{\iff}{\,\Leftrightarrow\,}
\newcommand{\argmax}{\mathop{\rm argmax}}
\newcommand{\pclass}{\mathscr{P}}
\newcommand{\tclass}{\mathscr{T}}
\newcommand{\idmat}{I}
\newcommand{\vecc}{\mathop{\mathrm{vec}}}
\newcommand{\lprod}{\mathop{\textstyle\prod}}
\newcommand{\normal}{\mathcal{N}}
\newcommand{\SSS}[1]{{\text{\small$#1$}}}
\newcommand{\parent}[1]{{\mathrm{pa}(#1)}}
\newcommand{\UE}{\mathrel{\relbar\joinrel\mkern-4.4mu\relbar\mkern-4.4mu\joinrel\relbar}}
\newcommand{\DE}{\mathrel{\relbar\joinrel\mkern-2mu\rightarrow}}
\newcommand{\LDE}{\mathrel{\leftarrow\joinrel\mkern-2mu\relbar}}
\begin{document}
\title
  {\bf Maximum Likelihood Estimation in Gaussian Chain Graph Models
  under the Alternative Markov Property} 
\author{Mathias Drton\\
    {\em Department of Statistics, The University of Chicago}
  \and Michael Eichler\\
  {\em Institute of Applied Mathematics, University of Heidelberg}
}
\date{}
\maketitle

\begin{abstract}
The AMP Markov property is a recently proposed alternative Markov
property for chain graphs. In the case of continuous variables with
a joint multivariate Gaussian distribution, it is the AMP rather
than the earlier introduced LWF Markov property that is coherent
with data-generation by natural block-recursive regressions. In this
paper, we show that maximum likelihood estimates in Gaussian AMP
chain graph models can be obtained by combining generalized least
squares and iterative proportional fitting to an iterative
algorithm.  In an appendix, we give useful convergence results for
iterative partial maximization algorithms that apply in particular
to the described algorithm.
\end{abstract}


\section{Introduction}

In graphical modelling, graphs are used to describe patterns of
conditional independence.  Undirected graphs encode the
conditional independences underlying Markov random fields, and acyclic
directed graphs encode the conditional independences underlying
Bayesian networks.  A generalization of both Markov random fields and
Bayesian networks is provided by chain graphs that were introduced
with the Markov/conditional independence interpretation described in
\cite{lauwerm:graphmod}, \cite{wermuth:1990} and
\cite{fryd:cgmp}; see also \citet[\S5.4.1]{lau:bk} and
\citet[\S7.2]{edwards:2000}.  Graphical models jargon refers to the models
induced by this Markov interpretation as LWF chain graph models.
Recently, however, \cite{amp:cgmp} have proposed an alternative Markov
property (AMP) for chain graphs \cite[see
also][]{levitz:2001,andersson:2004}.  In the case of continuous variables 
with a joint multivariate Gaussian = normal distribution, it is their
AMP rather than the LWF Markov property that is coherent with
data-generation by natural block-recursive regressions \citep[\S\S1
and 5]{amp:cgmp}.

Statistical inference for LWF chain graph models is well developed,
but this is not the case for AMP chain graph models.  This paper
considers maximum likelihood (ML) estimation in Gaussian AMP chain
graph models.  After reviewing these models in section
\ref{gaussian-amp}, we derive, in section \ref{max-likelihood}, the
likelihood equations and the Fisher-information.  Combining
generalized least squares and iterative proportional fitting, we
describe an iterative algorithm for solving the likelihood equations,
which yields consistent and asymptotically efficient estimates.  The
convergence properties of this algorithm can be derived from
convergence results for iterative partial maximization algorithms that
are given in the Appendix. An application to university graduation
data in section \ref{example} illustrates AMP chain graph modelling.
We conclude with the discussion in section \ref{sec:conc}.

\section{Gaussian AMP chain graph models}
\label{gaussian-amp}

Let $G=(V,E)$ be a mixed graph with finite vertex set $V$ and an edge
set $E$ that may contain two types of edges, namely directed ($u\DE
v$) and undirected ($u\UE v$) edges.  The graph $G$ is called a {\em chain
graph\/} if it does not contain any semi-directed cycles, that is, it
contains no path from $v$ to $v$ with at least one directed edge such
that all directed edges have the same orientation. The vertex set of a
chain graph can be partitioned into subsets $\tau\in\tclass$ such that
all edges within each subset $\tau$ are undirected and edges between
two different subsets $\tau\not=\tau'$ are directed.  In the
following, we assume that the partition $\tau\in\tclass$ is
maximal, that is, any two vertices in a subset $\tau$ are connected by
an undirected path. Then the subsets $\tau\in\tclass$ are unique and
called the {\em chain components\/} of the graph $G$; compare figure
\ref{fig:G4cycle} in section \ref{example}.

For a given chain graph $G$, we consider the class $\pclass(G)$ of
normal distributions $\normal(0,\Sigma)$ on $\rnum^V$ with positive
definite covariance matrix $\Sigma$ that satisfy the AMP Markov
property \citep[\S4]{amp:cgmp} with respect to $G$.
\citet[\S5]{amp:cgmp} described a parameterization of $\pclass(G)$
that associates one parameter with each vertex in $V$ and each edge in
$E$.  More precisely, let $\Omega=(\Omega_{uv})\in\rnum^{V\times V}$
be a positive definite matrix such that
for any distinct vertices $u$ and $v$
\begin{equation}
\label{constraints1}
u\UE v\notin E\follows\Omega_{uv}=0
\end{equation}
and let $B=(B_{uv})$ be an arbitrary matrix in $\rnum^{V\times V}$
such that for any vertices $u$ and $v$
\begin{equation}
\label{constraints2}
u\DE v\notin E\follows B_{vu}=0.
\end{equation}
For two such matrices $\Omega$ and $B$, we set
\begin{equation}\label{eq:sigma}
\Sigma(B,\Omega)=(\idmat_V-B)^{-1}\,\Omega^{-1}\,(\idmat_V-B')^{-1},
\end{equation}
where $\idmat_V\in\rnum^{V\times V}$ denotes the identity matrix.  A
normal distribution $\normal(0,\Sigma)$ with $\Sigma>0$ satisfies the
AMP Markov property if and only if there exist $B$ and $\Omega$ such
that \eqref{constraints1} and \eqref{constraints2} hold and
$\Sigma=\Sigma(B,\Omega)$.

For a vertex $v\in V$, let $\parent{v}=\{u\in V\mid u\DE v\in E\}$ be the
set of {\em parents\/} of $v$. Furthermore, we set
$\parent{\tau}=\cup_{v\in\tau}\parent{v}$.
Because of the nonexistence of semi-directed cycles, the joint
distribution of $X_V$ can be factorized as
\begin{equation}
\label{factorization}
f(x_V)=\lprod_{\tau\in\tclass} f(x_\tau\mid x_{\parent{\tau}}),
\quad x_V\in\rnum^V.
\end{equation}
For $\tau\in\tclass$, the conditional distribution
$f(x_\tau\mid x_{\parent{\tau}})$ is given by
\begin{equation}
\label{blockregression}
X_\tau\mid X_{\parent{\tau}}
\sim\normal\big(B_\tau\,X_{\parent{\tau}},\Omega_\tau^{-1}\big),
\end{equation}
where $B_\tau=(B_{uv})_{u\in\tau,v\in\parent{\tau}}$ and
$\Omega_\tau=(\Omega_{uv})_{u,v\in\tau}$ are submatrices of
$B$ and $\Omega$, respectively. The conditional distribution
corresponds to a {\em block-regression}, in which the block of variables
$X_\tau$ is regressed on the parents $X_{\parent{\tau}}$.

The parameter $(B_\tau,\Omega_\tau)$ can be rewritten in vectorized
form. Let $\beta_\tau=(B_{uv}\mid u\in\tau,v\in\parent{\tau})$ be the
vector of unconstrained elements in $B_\tau$.  Subsequently, we write
$B_\tau(\beta_\tau)$ for the matrix defined by $\beta_\tau$ and
\eqref{constraints2}.  Similarly, let $\omega_\tau$ be the vector of
elements of $\Omega_{uv}$ in $\Omega_\tau$ such that either $u=v$, or
$u<v$ and $u\UE v\in E$.
Furthermore, denote the dimension of $\beta_\tau$ and $\omega_\tau$
by $p_\tau$ and $q_\tau$, respectively. Then the parameter space
for the parameter $(\beta_\tau,\omega_\tau)$ is
\begin{equation}
\label{eq:parspace}
\Theta_\tau
=\big\{(\beta_\tau,\omega_\tau)\in\rnum^{p_\tau+q_\tau}\mid
\Omega_\tau(\omega_\tau)>0\big\},
\end{equation}
where $\Omega_\tau(\omega_\tau)\in\rnum^{\tau\times\tau}$ is the
matrix defined by $\omega_\tau$ and \eqref{constraints1}. It follows
from \eqref{factorization} and \eqref{blockregression} that
$\theta=(\beta_\tau,\omega_\tau)_{\tau\in\tclass}$ parameterizes
$\pclass(G)$.  Equation \eqref{eq:ident} below clarifies that $\theta$
is identifiable.  The parameter space of $\pclass(G)$ is the Cartesian
product $\Theta=\times_{\tau\in\tclass}\Theta_\tau$.  This
factorization of the parameter space together with the factorization
of the joint density implies that the ML estimator (MLE) of the joint
parameter $\theta$ can be obtained by computing, separately for every
$\tau\in\tclass$, the MLE of $(\beta_\tau,\omega_\tau)$ in the
block-regression \eqref{blockregression}. Furthermore, the Hessian of
the likelihood function of the model $\pclass(G)$ is block-diagonal
with one block for each one of the block-regressions indexed by
$\tau\in\tclass$.

\section{Maximum likelihood estimation}
\label{max-likelihood}

\subsection{Likelihood equations}
\label{sec:likeqn}

Let $X=(X_{v,i})_{v\in V,i\in N}\in\rnum^{V\times N}$ now be a data
matrix whose column vectors, indexed by the set $N$, are independent
and identically distributed according to some $P\in\pclass(G)$.
Since, merely for notational convenience, the distributions in
$\pclass(G)$ are assumed to be centered the sample covariance matrix
is defined as
\[ 
S= \frac{1}{n} X X',
\]
where $n=|N|$ is the sample size.  We assume that
\[ 
n\ge \max_{\tau\in\tclass}  \{ |\tau|+|\parent{\tau}| \}
\]
such that, with probability one, the submatrices $S_{\tau,\tau}$,
$S_{\tau,\parent{\tau}}$, and the matrix $S(\beta_\tau)$ defined below
are of full rank. This ensures that the MLE
exists in each one of the block-regressions.
Dividing by $n$ and ignoring the additive constant $-(|V|/2)\log(2\pi)$,
the log-likelihood function for the block-regression
\eqref{blockregression} is given by 
\begin{equation}
\label{likelihood}
\ell_n(\beta_\tau,\omega_\tau)
=\SSS{\frac{1}{2}}\,\log|\Omega_\tau(\omega_\tau)|
-\SSS{\frac{1}{2}}\,\trace\big[\Omega_\tau(\omega_\tau)\,S(\beta_\tau)\big],
\end{equation}
where
\[
\begin{split}
S(\beta_\tau)&=\SSS{\frac{1}{n}}\,
[X_\tau-B_\tau(\beta_\tau)\,X_{\parent{\tau}}]
[X_\tau-B_\tau(\beta_\tau)\,X_{\parent{\tau}}]'\\
&= S_{\tau,\tau} -B_\tau(\beta_\tau)S_{\parent{\tau},\tau}-
S_{\tau,\parent{\tau}}B_\tau(\beta_\tau)' +
B_\tau(\beta_\tau)S_{\parent{\tau},\parent{\tau}}B_\tau(\beta_\tau)'
\end{split}
\]
is the sample covariance matrix of the residuals in the
block-regression (\ref{blockregression}), and $X_A\in\rnum^{A\times
  N}$ denotes the submatrix of $X$ that comprises all rows with index
in $A$.

Let $P_\tau=\partial\vecc(B_\tau)/\partial\beta_\tau'$ and
$Q_\tau=\partial\vecc(\Omega_\tau)/\partial\omega_\tau'$. Both $P_\tau$
and $Q_\tau$ have entries in $\{0,1\}$ and satisfy
$\vecc(B_\tau)=P_\tau\beta_\tau$ and
$\vecc(\Omega_\tau)=Q_\tau\omega_\tau$, respectively.  Each column in  
$P_\tau$ has exactly one entry equal to one. A column in 
$Q_\tau$ has exactly one or exactly two entries equal to one depending
on whether the associated element in $\omega_\tau$ comes from the 
diagonal or the off-diagonal part of $\Omega_\tau$, respectively.
With these two matrices the likelihood equations obtained by
taking first derivatives with respect to $\beta_\tau$ and
$\omega_\tau$ can be written as
\begin{equation}
\label{betaeq}
P_\tau'\,\big[\vecc(\Omega_\tau\,S_{\tau,\parent{\tau}})-
(S_{\parent{\tau},\parent{\tau}}\otimes\Omega_\tau)\,P_\tau\,\beta_\tau\big]=0
\end{equation}
and 
\begin{equation}
\label{omegaeq}
Q_\tau'\,\vecc\big[\Omega_\tau^{-1}-S(\beta_\tau)\big]=0.
\end{equation}
Equation (\ref{omegaeq}) represents in a compact way the fact that the
covariance associated with an undirected edge in the AMP chain graph
is equal to its counterpart in $S(\beta_\tau)$, that is, it is equal
to the empirical covariance of residuals computed for fixed
$\beta_\tau$.  Thus, equation (\ref{omegaeq}) parallels the well-known
likelihood equations of undirected Gaussian graphical models.

\subsection{Two-step estimation}
\label{sec:twostep}

If every vertex in $\parent{\tau}$ is adjacent to all vertices in
$\tau$, then no constraints on $B_\tau$ are imposed and $P_\tau$
becomes an identity matrix. In this case the first set of equations
leads to the usual least squares estimator 
\begin{equation}
  \label{eq:ols}
\beta_\tau=
(S_{\parent{\tau},\parent{\tau}}^{-1}\otimes\idmat_\tau)
\vecc(S_{\tau,\parent{\tau}})
\iff   B_\tau=S_{\tau,\parent{\tau}}\,S_{\parent{\tau},\parent{\tau}}^{-1}.
\end{equation}
Thus the MLE of $(\beta_\tau,\omega_\tau)$ can be obtained by fitting
an undirected graph model to the residuals computed using the
regression coefficients estimates in \eqref{eq:ols}.  This can be done
using iterative proportional fitting
\cite[pp.~182--185]{speed:1986,whittaker:bk}, which generally will
terminate in finitely many steps only if the subgraph $G_\tau$ induced
by the chain component $\tau$ is decomposable.

In the case of general AMP chain graphs with constraints on $B_\tau$,
a similar two-step method can also be used for parameter estimation,
as described in \citet[\S7.5]{edwards:2000}:
\begin{enumerate}
\item[1.] estimate $\beta_\tau$ by least squares by regressing each
  $X_v$, $v\in\tau$, on its parents $X_\parent{v}$, 
\item[2.] estimate $\omega_\tau$ by fitting an undirected
graph model to the regression residuals.
\end{enumerate}
For general AMP chain graphs with restrictions on $B_\tau$, however,
the two equations (\ref{betaeq}) and (\ref{omegaeq}) for $\beta_\tau$
and $\omega_\tau$ cannot be solved separately and the MLE differs from
this two-step estimator.

\subsection{Algorithm for maximum likelihood estimation}
\label{sec:algo}

To compute the MLE, or rather a solution to the likelihood
equations, in the general case, we consider an iterative method based on 
alternately maximizing the likelihood with respect to
$\beta_\tau$ and $\omega_\tau$. Let $(\tilde\beta_\tau,\tilde\omega_\tau)$
be a consistent estimator. Then setting $\omega^{(1)}=\tilde\omega_\tau$
we define the sequence of estimators
\begin{equation}
\label{mle-seq1}
\hat\beta_\tau^{(k+1)}=\argmax_{\beta_\tau\in\rnum^{p_\tau}}
\ell_n(\beta_\tau,\hat\omega_\tau^{(k)})
\end{equation}
and
\begin{equation}
\label{mle-seq2}
\hat\omega_\tau^{(k+1)}=\argmax_{\Omega_\tau(\omega_\tau)>0}
\ell_n(\hat\beta_\tau^{(k+1)},\omega_\tau)
\end{equation}
for $k\geq 2$. Note that $\hat\beta_\tau^{(k+1)}$ can be computed in
an explicit formula as the solution to (\ref{betaeq}) with
$\Omega_\tau$ substituted by
$\Omega_\tau(\hat\omega_\tau^{(k)})$, which is
\begin{equation}
\label{mle-seq1-formula}
\hat\beta_\tau^{(k+1)}=
\big\{ P_\tau'\,[S_{\parent{\tau},\parent{\tau}}\otimes 
  \Omega_\tau(\hat\omega_\tau^{(k)})]\,P_\tau \big\}^{-1}
\big\{P_\tau'\,\vecc[\Omega_\tau(\hat\omega_\tau^{(k)})\,
  S_{\tau,\parent{\tau}}]\big\}. 
\end{equation}
Similarly, $\hat\omega_\tau^{(k+1)}$ can be computed as the solution
to (\ref{omegaeq}) with $\beta_\tau$ substituted by
$\hat\beta_\tau^{(k+1)}$.  The equations in (\ref{omegaeq}) then
correspond to the likelihood equations of an undirected graph model
for the undirected induced subgraph $G_\tau$ and the regression
residuals as data.  In other words the undirected graph model for
$G_\tau$ has to be fitted to the sample covariance matrix
$S(\hat\beta_\tau^{(k+1)})$, for which the iterative proportional
fitting algorithm can be used.

The convergence properties of the sequence
$(\hat\beta_\tau^{(k)},\hat\omega_\tau^{(k)})_{k\in\nnum}$ are
discussed in the appendix. In particular, it follows that the sequence
converges if there are only finitely many solutions to the likelihood
equations (\ref{betaeq}) and (\ref{omegaeq}).  Note that the
likelihood equations may indeed have multiple solutions; compare
\citet{drton:2004} and \cite{drton:2004b} who consider seemingly
unrelated regressions that are special cases of the block-regressions
encountered here.

\subsection{The Fisher-information}
\label{sec:fisher}

For $\tau\in\tclass$, the second derivatives of the log-likelihood
function are
\begin{align*}
\SSS{\frac{\partial^2\ell_n(\beta_\tau,\omega_\tau)}
{\partial\beta_\tau\partial\beta_\tau'}}
&=-P'_\tau\big[S_{\parent{\tau},\parent{\tau}}\otimes\Omega_{\tau}\big]P_\tau,\\
\SSS{\frac{\partial^2\ell_n(\beta_\tau,\omega_\tau)}
{\partial\omega_\tau\partial\omega_\tau'}}
&=-\SSS{\frac{1}{2}}\,Q'_\tau\big[\Omega_\tau^{-1}\otimes\Omega_\tau^{-1}\big]Q_\tau\\
\intertext{and}
\SSS{\frac{\partial^2\ell_n(\beta_\tau,\omega_\tau)}
{\partial\beta_\tau\partial\omega_\tau'}}
&=-P'_\tau\big[\big(S_{\parent{\tau},\tau}-S_{\parent{\tau},\parent{\tau}}
B_\tau(\beta_\tau)'\big)\otimes\idmat_\tau\big]Q_\tau.
\end{align*}
Let $\theta=(\theta_\tau)_{\tau\in\tclass}=
(\beta_\tau,\omega_\tau)_{\tau\in\tclass}\in\Theta$, and let $\Sigma$
be the associated covariance matrix given by equation
\eqref{eq:sigma}. Then the Fisher-information $\mathscr{I}(\theta)$
for the Gaussian AMP chain graph model $\pclass(G)$ is block-diagonal
and the $\tau\times\tau$-block is equal to
\begin{equation}
  \label{eq:fisher}
  \mathscr{I}(\theta)_{\tau,\tau} =
  \bigg(\begin{matrix}
    P'_\tau(\Sigma_{\parent{\tau},\parent{\tau}}\otimes\Omega_\tau)P_\tau &0\\
    0&\frac{1}{2}\,Q'_\tau\big(\Omega_\tau^{-1}\otimes
    \Omega_\tau^{-1}\big){Q_\tau} 
  \end{matrix}\bigg).
 \end{equation}

\subsection{Consistency and asymptotic normality}
\label{sec:consist}

In the following, let
$\hat\theta_{\tau,n}=(\hat\beta_{\tau,n},\hat\omega_{\tau,n})$ be the
limit of the sequence
$(\hat\beta_\tau^{(k)},\hat\omega_\tau^{(k)})_{k\in\nnum}$ for sample
size $n$ and let
$\hat\theta_n=(\hat\theta_{\tau,n})_{\tau\in\tclass}$. Should such a
limit not exist then choose $\hat\theta_{\tau,n}$ as an arbitrary
accumulation point.  In either situation, all $\hat\theta_{\tau,n}$
are roots to the likelihood equations (\ref{betaeq}) and
(\ref{omegaeq}).  This, together with the fact that Gaussian AMP chain
graph models form curved exponential families (theorem \ref{thm:cef}),
leads to the asymptotic normality stated in theorem \ref{thm:asy}.

\begin{theorem}
  \label{thm:cef}
  The Gaussian AMP chain graph model $\pclass(G)$ is a curved
  exponential family.
\end{theorem}
{\em Proof.}  The model $\pclass(G)$ is a subfamily of the regular
exponential family of centered multivariate normal distributions with
arbitrary positive definite covariance matrix.  The parameter space
$\Theta=\times_{\tau\in\tclass} \Theta_\tau$ of $\pclass(G)$ is an
open set in a Euclidian space and in particular a smooth manifold.  For
$\theta=(\beta_\tau,\omega_\tau)_{\tau\in\tclass}\in\Theta$, let
$B(\theta)$ be the matrix that is zero except for its
$\tau\times\parent{\tau}$-submatrices, $\tau\in\tclass$, which are
equal to $B_\tau(\beta_\tau)$, and let similarly $\Omega(\theta)$ be
the block-diagonal matrix with blocks $\Omega_\tau(\omega_\tau)$,
$\tau\in\tclass$.  By equation (\ref{eq:sigma}), the mapping
\[
\psi:\theta\mapsto \Sigma(\theta)^{-1} =
[\idmat_V-B(\theta)']\,\Omega(\theta)\,[\idmat_V-B(\theta)] 
\]
maps the parameter
$\theta=(\beta_\tau,\omega_\tau)_{\tau\in\tclass}\in\Theta$ in the
parameter space of $\pclass(G)$ to $\Sigma(\theta)^{-1}$ with
$\normal(0,\Sigma(\theta))\in\pclass(G)$.  The inverse map of $\psi$
is determined by the fact that
\begin{equation}
\label{eq:ident}
B(\theta)_{v,\parent{v}}=
\Sigma(\theta)_{v,\parent{v}}
[\Sigma(\theta)_{\parent{v},\parent{v}}]^{-1},\qquad v\in V;
\end{equation}
compare \citet[Theorem~8.7]{richardson:2002}.  It is now apparent that
the mapping $\psi$ is a diffeomorphism.
Therefore, $\psi(\Theta)$ is a smooth manifold, which means that
$\pclass(G)$ forms a curved exponential family \cite[Definition~2.3.1,
4.2.1]{kassvos}.
\bigskip

\begin{theorem}
  \label{thm:asy}
  Let $\theta=(\theta_\tau)_{\tau\in\tclass}$,
  $\theta_\tau=(\beta_\tau,\omega_\tau)$, be the true parameter. Then
  $\hat\theta_n\to\theta$ in probability, the estimates
  $\hat\theta_{\tau,n}$, $\tau\in\tclass$, are asymptotically
  independent, and for each $\tau\in\tclass$,
\[
\sqrt{n}\begin{pmatrix}
\hat{\beta}_{\tau,n}-\beta_\tau\\
\hat{\omega}_{\tau,n}-\omega_\tau
\end{pmatrix}
\dconv\normal\big( 0, [\mathscr{I}(\theta)_{\tau,\tau}]^{-1}\big)
\]
with $\mathscr{I}(\theta)_{\tau,\tau}$ given in (\ref{eq:fisher}).
\end{theorem}
{\em Proof.}  The estimators $\hat\theta_n$ are roots to the likelihood
equations, computed in iterations starting at consistent estimates.
Theorems 2.4.1, 2.6.1, 2.6.7, and 2.6.12 (see also Corollaries 2.4.2 and
2.6.2) in \citet{kassvos} imply that in one-parameter curved exponential
families such roots to the likelihood equations are consistent and
asymptotically normal with asymptotic variance equal to the inverse of the
Fisher-information.  As indicated before the statement of Theorem 4.2.4 in
\citet{kassvos}, these results extend to multi-parameter families, which
yields our claim.
\bigskip

\section{Example: University graduation rates}
\label{example}

We illustrate our maximum likelihood procedure using the data in
\cite{druzdzel:1999}, which stem from a study
for college ranking carried out in 1993.  Based on $n=159$
universities, \citet[Table 3]{druzdzel:1999} state a correlation
matrix for eight variables that are
\begin{center}
\begin{tabular}{rl}
\em spend & average spending per student,\\
\em strat & student-teacher ratio,\\
\em salar & faculty salary,\\
\em rejr & rejection rate,\\
\em pacc & percentage of admitted students who accept university's offer,\\
\em tstsc & average test scores of incoming students,\\
\em top10 & class standing of incoming freshmen, and\\
\em apgra & average percentage of graduation.
\end{tabular}
\end{center}

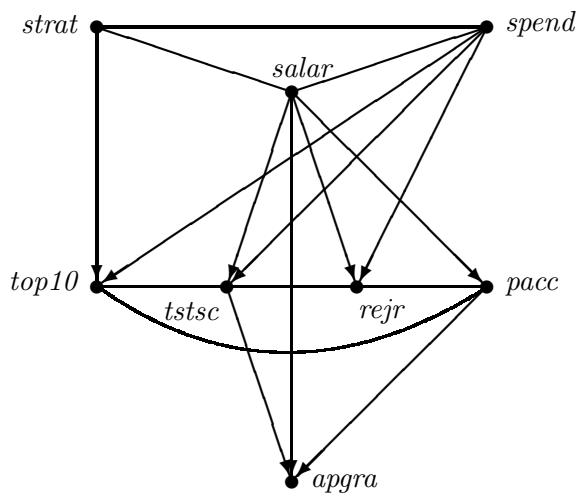
\begin{figure}[t]
  \centering
  \setlength{\unitlength}{0.55pt}
  \begin{picture}(367,429)(0,80)
    \thicklines
    \put(323.2,407.0){\vector(-3,-2){   261.5}}
    \put(323.7,406.4){\vector(-1,-1){   173.1}}
    \put(324.7,405.7){\vector(-1,-2){    85.8}}
    \put(191.2,360.8){\vector(-1,-3){    42.2}}
    \put(193.7,360.8){\vector(1,-3){    42.2}}
    \put(195.3,361.7){\vector(1,-1){   128.4}}
    \put(192.4,360.6){\vector(0,-1){   260.2}}
    \put( 58.3,405.3){\vector(0,-1){   170.8}}
    \put(149.0,226.7){\vector(1,-3){    42.2}}
    \put(323.7,227.7){\vector(-1,-1){   128.4}}
    \put(322.7,408.0){\line(-3,-1){   126.5}}
    \put(322.5,409.3){\line(-1,0){   260.2}}
    \put(188.6,365.8){\line(-3,1){   126.5}}
    \put( 62.3,230.5){\line(1,0){    81.4}}
    \put(151.7,230.5){\line(1,0){    81.4}}
    \put(241.1,230.5){\line(1,0){    81.4}}
    \qbezier(58.3,230.5)(184.65,140)(326.5,230.5)
    \put(326.5,409.3){\circle*{9}}
    \put(339,403){\makebox(0,0)[l,b]{\em spend}}
    \put(192.4,364.6){\circle*{9}}
    \put(178,375){\makebox(0,0)[l,b]{\em salar}}
    \put( 58.3,409.3){\circle*{9}}
    \put(45,405){\makebox(0,0)[r,b]{\em strat}}
    \put( 58.3,230.5){\circle*{9}}
    \put(45,224.2){\makebox(0,0)[r,b]{\em top10}}
    \put(147.7,230.5){\circle*{9}}
    \put(142,220){\makebox(0,0)[r,t]{\em tstsc}}
    \put(237.1,230.5){\circle*{9}}
    \put(238,222){\makebox(0,0)[l,t]{\em rejr}}
    \put(326.5,230.5){\circle*{9}}
    \put(339,224.2){\makebox(0,0)[l,b]{\em pacc}}
    \put(192.4, 96.4){\circle*{9}}
    \put(205,89){\makebox(0,0)[l,b]{\em apgra}}
  \end{picture}
  \bigskip
  \normalsize
  \caption{Chain graph with the three chain components $\{\mbox{\em spend,
      strat, salar}\}$,  $\{\mbox{\em top10, tstsc, rejr, pacc}\}$, and
    $\{\mbox{\em agpra}\}$.}  
\label{fig:G4cycle}
\end{figure}

Figure \ref{fig:G4cycle} shows a chain graph for these variables.
This graph has the chain components $\tau_1=\{\mbox{\em spend, strat,
  salar}\}$, $\tau_2=\{\mbox{\em top10, tstsc, rejr, pacc}\}$, and
$\tau_3=\{\mbox{\em apgra}\}$.  It was selected via the SIN model
selection procedure described in
\cite{drtonperlman:2004,drtonperlman:2004b}.  More precisely, we used
SIN model selection with simultaneous significance level $0.15$ fixing
the chain components $\tau_1$, $\tau_2$, and $\tau_3$ {\em a priori\/}
in the temporal order $\tau_1<\tau_2<\tau_3$.  In the resulting AMP
chain graph we deleted the undirected edge between {\em top10} and
{\em rejr}, and introduced the undirected edge between {\em top10} and
{\em pacc}, creating a non-decomposable chain component $\tau_2$.
Furthermore, we deleted the edge between {\em salar} and {\em top10}
to create the edge constellation
\[
\mbox{\em salar}\DE\mbox{\em top10}\UE \mbox{\em tstsc}\LDE\mbox{\em strat}.
\]
The induced subgraph over the four vertices {\em salar}, {\em strat},
{\em top10} and {\em tstsc}, which also contains the edge $\mbox{\em
  salar}\UE\mbox{\em strat}$, forms what is called a 2-biflag by
\cite{amp:cgmp}; compare their figure~5(d).  Therefore, by theorem 4
in \cite{amp:cgmp}, the AMP and LWF Markov properties differ for the
graph in figure \ref{fig:G4cycle}.

The block-regression for $\tau_1$ is trivial as $\parent{\tau_1}=\emptyset$
and the undirected induced subgraph $G_{\tau_1}$ is complete, and thus the
MLE of $\Omega_{\tau_1}$ is simply the inverse of $S_{\tau_1,\tau_1}$.  The
block-regression for $\tau_3$ is also simple as $\tau_3$ contains only a
single vertex.  In this case, the MLE of $\beta_{\tau_3}$ and
$\omega_{\tau_3}$ can be computed by regressing the single variable in
$\tau_3$, here the variable {\em apgra}, on all its parents, here the
variables {\em pacc}, {\em salar}, and {\em tstsc}.  The vector of least
squares estimates of the regression coefficients is the MLE of
$\beta_{\tau_3}$ and the inverse of the estimated conditional variance is
the MLE of $\omega_{\tau_3}$.

The remaining block-regression for $\tau_2$ is non-trivial.  We apply
the ML estimation algorithm described in section \ref{sec:algo},
starting with the identity matrix as initial estimate of
$\Omega_{\tau_2}$ and iterating until convergence to find the
estimates stated in the columns labelled ``MLE'' in table
\ref{tab:estim}. Note that we cannot guarantee that these estimates
constitute the global maximum of the likelihood function.  However,
using these estimates to evaluate the deviance of the AMP chain graph
model yields a value of 16.89, which compared to 11 degrees of freedom
indicates a reasonable fit.

\begin{table}[htb]
  \centering
  \caption{\em MLEs, their standard errors computed from the
    Fisher information matrix, and the two-step
    estimates for the block-regression for chain-component $\tau_2$.}
  \label{tab:estim}
  \vspace{0.9cm}
  \begin{tabular}{crrrr@{\hspace{1.5cm}}crrr} 
    \hline
    Parameter & \multicolumn{1}{c}{MLE} & \multicolumn{1}{c}{SE} &
    \multicolumn{1}{c}{2-step} &&
    Parameter & \multicolumn{1}{c}{MLE} & \multicolumn{1}{c}{SE} &
    \multicolumn{1}{c}{2-step}\\ 
    \hline
    $\beta_{\mathit{pacc}\leftarrow\mathit{salar}}$  & -0.53 & 0.07 & -0.52 &&
    $\omega_{\mathit{pacc}}$  & 1.46 & 0.16 & 1.46\\
    $\beta_{\mathit{rejr}\leftarrow\mathit{salar}}$  &  0.26 & 0.09 & 0.30&&
    $\omega_{\mathit{rejr}}$  & 1.64 & 0.18 & 1.64\\
    $\beta_{\mathit{rejr}\leftarrow\mathit{spend}}$  &  0.30 & 0.09 & 0.27&&
    $\omega_{\mathit{top10}}$ & 2.99 & 0.33 & 2.92\\
    $\beta_{\mathit{top10}\leftarrow\mathit{spend}}$ &  0.98 & 0.08 & 0.99&&
    $\omega_{\mathit{tstsc}}$ & 3.39 & 0.37 & 3.34\\
    $\beta_{\mathit{top10}\leftarrow\mathit{strat}}$ &  0.44 & 0.07 & 0.45&&
    $\omega_{\mathit{pacc},\mathit{rejr}}$   & -0.33 & 0.12 & -0.33\\
    $\beta_{\mathit{tstsc}\leftarrow\mathit{salar}}$ &  0.26 & 0.06 & 0.36&&
    $\omega_{\mathit{pacc},\mathit{top10}}$  & -0.16 & 0.14 & -0.16\\
    $\beta_{\mathit{tstsc}\leftarrow\mathit{spend}}$ &  0.49 & 0.07 & 0.43&&
    $\omega_{\mathit{rejr},\mathit{tstsc}}$  & -0.65 & 0.16 & -0.65\\
    & & & &&
    $\omega_{\mathit{top10},\mathit{tstsc}}$ & -1.76 & 0.28 & -1.69\\
    \hline
  \end{tabular}
\end{table}

Table \ref{tab:estim} also states the two-step estimates obtained as
described in section \ref{sec:twostep}.  These estimates coincide with
the estimates after two steps of the ML estimation algorithm, provided
the algorithm is started at a diagonal matrix.  The two steps of the
ML estimation algorithm consist of one step estimating $\beta_\tau$
assuming a diagonal matrix $\Omega_\tau$, i.e.~assuming independence
of all variables in the chain component $\tau$, and one step
estimating $\omega_\tau$ using the newly found estimate of
$\beta_\tau$.  The two-step estimates are fairly close to the MLEs,
all differences being clearly smaller than two standard errors.  The
deviance based on the two-step estimates would be 19.18.
Interestingly, the two-step estimates and the MLEs for the variance
parameters $\omega_\tau$ are identical in two digits of precision with
the exception of the conditional variances $\omega_{\mathit{top10}}$,
$\omega_{\mathit{tstsc}}$ and the inverse covariance
$\omega_{\mathit{top10},\mathit{tstsc}}$ that all involve the
variables {\em top10} and {\em tstsc} that are part of the biflag.

\section{Discussion}
\label{sec:conc}

The likelihood function of a Gaussian AMP chain graph model can be
factored into the product of conditional likelihood functions.  Each
chain component of the graph gives rise to one factor in this
factorization.  The iterative algorithm we proposed for ML estimation
in Gaussian AMP chain graph models takes advantage of this fact and
treats each chain component separately.  For a given chain component,
the algorithm alternates between estimating regression coefficients
while fixing a covariance matrix and estimating the (restricted)
covariance matrix while fixing regression coefficients.  To perform
the former task of estimating regression coefficients we use a
generalized least squares formula, whereas the iterative proportional
fitting algorithm is used to perform the latter task of estimating a
covariance matrix.

The algorithm calls upon repeated runs of iterative proportional
fitting in order to fit the block-regression model associated with a
given chain component.  This is in contrast to the case of LWF chain
graph models, for which the ML estimates of the parameters associated
with a chain component can be computed by running iterative
proportional fitting only once \citep[\S 5.4.1, Proposition~6.33]{lau:bk}.
However, the undirected graph on which iterative proportional fitting
is run must be derived from the original LWF chain graph in a process
called moralization.  In general, this derived undirected graph
contains also vertices outside the considered chain component and may
feature larger cliques than the undirected subgraph induced by the
chain component, on which iterative proportional fitting is run when
fitting AMP chain graph models.

The developed methodology for ML fitting of AMP chain graph models
permits in particular to compare two models based on different chain
graphs via likelihood ratio tests and information criteria.  However,
one may also be interested in testing parameter equality in a given
model.  If parameters are set equal in a curved exponential family,
then the resulting submodel is again a curved exponential family.
Therefore, the ML estimates in the submodel are asymptotically normal,
and the problem of testing parameter equality can be addressed by a
likelihood ratio test.  For the computation of ML estimates in such
submodels, the algorithm we proposed for fitting AMP chain graph
models needs to be extended to incorporate equality constraints
amongst subsets of the parameters.  If parameter equality occurs
between regression coefficients that appear in the same matrix
$B_\tau$, then the generalized least squares step of the algorithm can
easily be adapted to deal with this new situation.  The required
changes consist solely of removing all but one of the identical
entries of the vector $\beta_\tau$ and altering the matrix $P_\tau$
accordingly.  With these changes, formula \eqref{mle-seq1} still
applies.  If parameter equality occurs between entries of the matrix
$\Omega_\tau$ then the iterative proportional fitting step of the
algorithm has to be adapted.  This can be done as described in
\cite{hojsgaard:2005} who treat parameter equality in undirected
graphical models.  Finally if parameter equality occurs between
parameters appearing in different matrices $B_\tau$ and
$B_{\bar\tau}$, or $\Omega_\tau$ and $\Omega_{\bar\tau}$, then the
block-regressions can no longer be treated separately. In this case
the extension of the presented algorithm requires additional work.

\section*{Appendix: Iterative partial maximization}
The algorithm for ML estimation proposed in this
paper is an iterative partial maximization algorithm in the sense of
\citet[Appendix A.4]{lau:bk}. Partial maximization refers to a
maximization of the likelihood function over a section in the
parameter space.  In an  iterative partial maximization algorithm, one
repeatedly performs a sequence of partial maximizations. 
In this appendix, we generalize the convergence
results in  \citet[Appendix A.4]{lau:bk} by not assuming the existence
of a unique local (and global) maximum of the likelihood function.

Let $L:\theta\to\rnum$ be a differentiable real-valued function on an open
set $\Theta\subseteq\rnum^d$. In the context of ML estimation, $L$
constitutes the (log-)likelihood function and $\Theta$ is the parameter
space of a statistical model. Assume that there exists $\theta_0$ such that
$\Theta_0=\{\theta\in\Theta\mid L(\theta)\geq L(\theta_0)\}$ is compact.
Then $L$ has a (not necessarily unique) global maximum in $\Theta_0$.  For
functions $g_i:\Theta\to\rnum^{d_i}$, $i=1,\ldots,k$ and
$\theta^*\in\Theta$, we define sections $\Theta_i(\theta^*)$ in $\Theta$ by
\[
\Theta_i(\theta^*)
=\{\theta\in\Theta\mid g_i(\theta)=g_i(\theta^*)\}.
\]
We assume that the maximum of $L$ over the section
$\Theta_i(\theta^*)$ is uniquely attained for all $\theta^*\in\Theta$ and
$i=1,\ldots,k$ and that the associated mapping
\[
T_i(\theta^*)=\argmax_{\theta\in\Theta_i(\theta^*)}L(\theta)
\]
from $\Theta$ into itself is continuous for all $i=1,\ldots,k$.
Moreover, we assume that if $\theta^*$ maximizes $L$ over all sections
$\Theta_i(\theta^*)$, and consequently satisfies
$\theta^*=T_i(\theta^*)$ for $i=1,\ldots,k$, then $\theta^*$ solves
the likelihood equations
\begin{equation}
\label{L equations}
\frac{\partial L(\theta)}{\partial\theta}\bigg|_{\theta=\theta^*}=0.
\end{equation}

\newcommand{\acset}{\mathcal{A}_\infty}

Let $\theta_0\in\Theta$ be a starting value such that
$\Theta_0$ is compact and define
\[
\theta_{n+1}=S(\theta_n)=T_k\cdots T_1(\theta_n),\qquad n\geq 0.
\]
By definition of $\Theta_0$, we have $\theta_n\in\Theta_0$ for all
$n\geq 0$. Let $\acset$ be the set of accumulation points of
the sequence $(\theta_n)_{n\in\nnum}$. Since $\Theta_0$ is compact,
we have $\acset\subseteq\Theta_0$. The following results
discuss the properties of $\acset$. It is a special case of the convergence
theorem in \citet[Chapter~4]{zangwill:1969}.

\begin{proposition}
\label{accum-prop}
The sequence $\big(L(\theta_n)\big)_{n\in\nnum}$ of values of the
likelihood function converges to a limit $\ell_\infty\in\rnum$.
Furthermore, if $\alpha\in\acset$ then $L(\alpha)=\ell_\infty$ and
$\alpha$ satisfies \eqref{L equations}.
\end{proposition}
{\em Proof.}
Since the sequence $(L(\theta_n))_{n\in\nnum}$ is monotonously increasing
and bounded, it converges to a limit $\ell_\infty$. By continuity
of $L$, this also implies $L(\alpha)=\ell_\infty$ for all $\alpha\in\acset$.

Next, since $S=T_k\cdots T_1$ is continuous, $S(\acset)$ is the set
of accumulation points of $(S(\theta_n))=(\theta_{n+1})$.
Consequently $S(\acset)=\acset$ and $L(S(\alpha))=\ell_\infty$ for
all $\alpha\in\acset$. By definition of $T_i$, we now obtain for
arbitrary $\alpha\in\acset$,
\[
\ell_\infty = L\big(T_k\cdots T_1(\alpha)\big)
\geq L\big(T_{k-1}\cdots T_1(\alpha)\big)
\geq L\big(T_1(\alpha)\big)
\geq L\left(\alpha\right)=\ell_\infty,
\]
which implies $T_i(\alpha)=\alpha$ for all $i=1,\ldots,k$
because of uniqueness of the maximum over $\Theta_i(\alpha)$.
Thus $\alpha$ maximizes $L$ over all sections and hence
 satisfies equations \eqref{L equations} by assumption.
\bigskip

For the next theorem, recall that a compact set is said to be connected
if it cannot be partitioned into two nonempty compact sets 
\cite[see also][Theorem~28.1]{ostrowski:1966}. 

\begin{theorem}
\label{accum-theorem}
$\acset$ is a compact and connected subset of $\Theta_0$.
\end{theorem}
{\em Proof.}
Since $\acset$ is a subset of a compact set, it suffices
to show that $\acset$ is closed. Let $\alpha^*\in\overline\acset$.
Then for any $\veps>0$ there exists $\alpha\in\acset$
such that $\alpha\in B_\veps(\alpha^*)$. Similarly, since
$\alpha$ is an accumulation point of $(\theta_n)$, there exists
for every $\delta>0$ some $n_\delta\in\nnum$
such that $\theta_{n_\delta}\in B_\delta(\alpha)$.
Since $B_\veps(\alpha^*)$ is open, we can choose $\delta$ small enough
such that $B_\delta(\alpha)\subseteq B_\veps(\alpha^*)$, which implies
$\theta_{n_\delta}\in B_\veps(\alpha^*)$. Since $\veps$ was arbitrary,
$\alpha^*\in \acset$, which establishes the closedness of $\acset$.

Next, let $B_\veps(\acset)=\{\theta\in\Theta \mid d(\theta,\acset)<\veps\}$
where $d(A,B)$ 
is the distance between two subsets $A$ and $B$ in $\rnum^d$.
Then for every $\veps>0$ there exists $n_\veps\in\nnum$ such that
$\theta_n\in B_\veps(\acset)$ for all $n\geq n_\veps$.

Now suppose that $\acset$ can be partitioned into
two compact sets $A$ and $B$. Then
$d(A,B)>0$ and we set $\delta=d(A,B)/2$.
Furthermore, because of uniform continuity of $S$ on $\Theta_0$,
for all $\delta>0$ there exists $\veps'>0$ such that for all
$\alpha\in\acset$, $\theta_n\in B_{\veps'}(\alpha)$ implies
$\theta_{n+1}=S(\theta_n)\in B_\delta(\alpha)$. 

Then if $n>n_\veps$ and $\theta_n\in B_\veps(A)$, we have
\[
\theta_{n+1}\in B_\delta(A)\cap B_\veps(\acset)=B_\veps(A),
\]
since $d(A,B)>\delta$. Thus $\theta_n\in B_\veps(A)$ for all
$n>n_\veps$ and hence $B=\varnothing$ which concludes the proof.
\bigskip

\begin{corollary}
\label{accum-corr1}
If $\acset$ is finite, then $\acset=\{\theta^*\}$
for some $\theta^*\in\Theta_0$ and the sequence
$(\theta_n)_{\theta\in\nnum}$ converges to $\theta^*$.
\end{corollary}
{\em Proof.}  Any connected finite set must be a singleton.
\bigskip

\begin{corollary}
\label{accum-corr2}
If the likelihood equations \eqref{L equations} have only finitely
many solutions that lie on the same contour of the likelihood
function $L$, then the sequence $(\theta_n)_{n\in\nnum}$
converges to one solution $\theta^*$.
\end{corollary}
{\em Proof.}
This follows from Proposition \ref{accum-prop} and
Corollary \ref{accum-corr1}.

\begin{small}

\end{small}

\bigskip
\bigskip
\noindent
Mathias Drton,
Department of Statistics,
The University of Chicago, 
5734 S. University Avenue,
Chicago IL, 60637, U.S.A.\\
E-mail: drton@galton.uchicago.edu


\begin{thebibliography}{}

\bibitem[\protect\citeauthoryear{Andersson, Madigan, \& Perlman}{Andersson {\it
  et~al.}}{2001}]{amp:cgmp}
Andersson, S.~A., Madigan, D. \& Perlman, M.~D. (2001).
\newblock Alternative {M}arkov properties for chain graphs.
\newblock {\em Scand. J. Statist.\/}~{\bf 28}, 33--85.

\bibitem[\protect\citeauthoryear{Andersson \& Perlman}{Andersson \&
  Perlman}{2004}]{andersson:2004}
Andersson, S.~A. \& Perlman, M.~D. (2004).
\newblock Characterizing {M}arkov equivalence classes for {AMP} chain graph
  models.
\newblock Technical Report 453, Department of Statistics, University of
  Washington.\\
\newblock Available at http://www.stat.washington.edu/www/research/reports/.

\bibitem[\protect\citeauthoryear{Drton}{Drton}{2005}]{drton:2004b}
Drton, M. (2005).
\newblock Computing all roots of the likelihood equations of seemingly
  unrelated regressions.
\newblock {\em J. Symbolic Comput.\/}, in press.

\bibitem[\protect\citeauthoryear{Drton \& Perlman}{Drton \&
  Perlman}{2004a}]{drtonperlman:2004}
Drton, M. \& Perlman, M.~D. (2004a).
\newblock Model selection for {G}aussian concentration graphs.
\newblock {\em Biometrika\/}~{\bf 91}, 591--602.

\bibitem[\protect\citeauthoryear{Drton \& Perlman}{Drton \&
  Perlman}{2004b}]{drtonperlman:2004b}
Drton, M. \& Perlman, M.~D. (2004b).
\newblock A {SIN}ful approach to {G}aussian graphical model selection.
\newblock Technical Report 457, University of Washington.\\
\newblock Available at http://www.stat.washington.edu/www/research/reports/.

\bibitem[\protect\citeauthoryear{Drton \& Richardson}{Drton \&
  Richardson}{2004}]{drton:2004}
Drton, M. \& Richardson, T.~S. (2004).
\newblock Multimodality of the likelihood in the bivariate seemingly unrelated
  regressions model.
\newblock {\em Biometrika\/}~{\bf 91}, 383--392.

\bibitem[\protect\citeauthoryear{Druzdzel \& Glymour}{Druzdzel \&
  Glymour}{1999}]{druzdzel:1999}
Druzdzel, M.~J. \& Glymour, C. (1999).
\newblock Causal inferences from databases: Why universities lose students.
\newblock In C.~Glymour \& G.~F. Cooper (Eds.), {\em Computation, {C}ausation,
  and {D}iscovery}, chapter~19, pp.\  521--539. AAAI Press, Menlo Park, CA.

\bibitem[\protect\citeauthoryear{Edwards}{Edwards}{2000}]{edwards:2000}
Edwards, D.~M. (2000).
\newblock {\em Introduction to graphical modelling\/} (Second ed.).
\newblock Springer-Verlag, New York.

\bibitem[\protect\citeauthoryear{Frydenberg}{Frydenberg}{1990}]{fryd:cgmp}
Frydenberg, M. (1990).
\newblock The chain graph {M}arkov property.
\newblock {\em Scand. J. Statist.\/}~{\bf 17}, 333--353.

\bibitem[\protect\citeauthoryear{H{\o}jsgaard \& Lauritzen}{H{\o}jsgaard \&
  Lauritzen}{2005}]{hojsgaard:2005}
H{\o}jsgaard, S. \& Lauritzen, S. (2005).
\newblock Restricted concentration models -- {G}aussian models with
  concentration parameters restricted to being equal.
\newblock In R.~G. Cowell \& Z.~Ghahramani (Eds.), {\em Proceedings of the
  Tenth International Workshop on Artificial Intelligence and Statistics}, pp.\
   152--157. Society for Artificial Intelligence and Statistics.\\
\newblock Available at http://www.gatsby.ucl.ac.uk/aistats/.

\bibitem[\protect\citeauthoryear{Kass \& Vos}{Kass \& Vos}{1997}]{kassvos}
Kass, R.~E. \& Vos, P.~W. (1997).
\newblock {\em Geometrical foundations of asymptotic inference}.
\newblock Wiley, New York.

\bibitem[\protect\citeauthoryear{Lauritzen}{Lauritzen}{1996}]{lau:bk}
Lauritzen, S.~L. (1996).
\newblock {\em Graphical models}.
\newblock Clarendon Press, Oxford, UK.

\bibitem[\protect\citeauthoryear{Lauritzen \& Wermuth}{Lauritzen \&
  Wermuth}{1989}]{lauwerm:graphmod}
Lauritzen, S.~L. \& Wermuth, N. (1989).
\newblock Graphical models for association between variables, some of which are
  qualitative and some quantitative.
\newblock {\em Ann. Statist.\/}~{\bf 17}, 31--57.

\bibitem[\protect\citeauthoryear{Levitz, Perlman, \& Madigan}{Levitz {\it
  et~al.}}{2001}]{levitz:2001}
Levitz, M., Perlman, M.~D. \& Madigan, D. (2001).
\newblock Separation and completeness properties for {AMP} chain graph {M}arkov
  models.
\newblock {\em Ann. Statist.\/}~{\bf 29}, 1751--1784.

\bibitem[\protect\citeauthoryear{Ostrowski}{Ostrowski}{1966}]{ostrowski:1966}
Ostrowski, A.~M. (1966).
\newblock {\em Solution of equations and systems of equations}.
\newblock Academic Press, New York.

\bibitem[\protect\citeauthoryear{Richardson \& Spirtes}{Richardson \&
  Spirtes}{2002}]{richardson:2002}
Richardson, T.~S. \& Spirtes, P. (2002).
\newblock Ancestral graph {M}arkov models.
\newblock {\em Ann.~Statist.\/}~{\bf 30}, 962--1030.

\bibitem[\protect\citeauthoryear{Speed \& Kiiveri}{Speed \&
  Kiiveri}{1986}]{speed:1986}
Speed, T.~P. \& Kiiveri, H.~T. (1986).
\newblock Gaussian {M}arkov distributions over finite graphs.
\newblock {\em Ann. Statist.\/}~{\bf 14}, 138--150.

\bibitem[\protect\citeauthoryear{Wermuth \& Lauritzen}{Wermuth \&
  Lauritzen}{1990}]{wermuth:1990}
Wermuth, N. \& Lauritzen, S.~L. (1990).
\newblock On substantive research hypotheses, conditional independence graphs
  and graphical chain models.
\newblock {\em J. Roy. Statist. Soc. Ser. B\/}~{\bf 52}, 21--50, 51--72.

\bibitem[\protect\citeauthoryear{Whittaker}{Whittaker}{1990}]{whittaker:bk}
Whittaker, J. (1990).
\newblock {\em Graphical models in applied multivariate statistics}.
\newblock Wiley, Chichester.

\bibitem[\protect\citeauthoryear{Zangwill}{Zangwill}{1969}]{zangwill:1969}
Zangwill, W.~I. (1969).
\newblock {\em Nonlinear programming: A unified approach}.
\newblock Prentice-Hall Inc., Englewood Cliffs, NJ.

\end{thebibliography}
\end{document}